\documentclass[12pt]{article}
\input{epsf.sty}
\usepackage{epsfig}
\usepackage{amsmath}
\usepackage{latexsym}
\usepackage{amssymb}
\usepackage{graphicx}
\usepackage{dcolumn}
\usepackage{bm}
\usepackage{times}

\newtheorem{theorem}{Theorem}
\newtheorem{lemma}{Lemma}

\newtheorem{definition}{Definition}

\newtheorem{example}{Example}
\newcommand{\onetom}{1,\ldots,m}
\newcommand{\oneton}{1,\ldots,n}

\newtheorem{myprop}{Proposition}

\def\sqr#1#2{{\vcenter{\vbox{\hrule height.#2pt \hbox{\vrule width.
#2pt height#1pt \kern#1pt \vrule} \hrule height.#2pt}}}}

\begin{document}
\begin{center}

{\large\bf Continuation of solutions of coupled dynamical
systems}\footnote{This work is supported by National Science
Foundation of China 60374018, 60574044. }
\\[0.2in]
Tianping Chen\footnote{Email:tchen@fudan.edu.cn}, Wei
Wu\footnote{051018023@fudan.edu.cn. These authors are with Lab. of
Nonlinear Mathematics Science, Institute of Mathematics, Fudan
University, Shanghai, 200433, P.R.China.\\
\indent ~~Corresponding author: Tianping Chen.
Email:tchen@fudan.edu.cn}

\end{center}

\begin{abstract}
Recently, the synchronization of coupled dynamical systems has been
widely studied. Synchronization is referred to as a process wherein
two (or many) dynamical systems are adjusted to a common behavior as
time goes to infinity, due to coupling or forcing. Therefore, before
discussing synchronization, a basic problem on continuation of the
solution must be solved: For given initial conditions, can the
solution of coupled dynamical systems be extended to the infinite
interval $[0,+\infty)$? In this paper, we propose a general model of
coupled dynamical systems, which includes previously studied systems
as special cases, and prove that under the assumption of QUAD, the
solution of the general model exists on $[0,+\infty)$.
\end{abstract}

Coupled dynamical systems, Synchronization, Existence, Uniqueness,
Continuation.

\newcounter{mycounter}
\section{Introduction}

In past years, collective behaviors of coupled dynamical systems
have been widely studied. In particular, synchronization in networks
of coupled dynamical systems, as one of the simplest and most
striking behaviors, has attracted increasing attention in
mathematical and physical literatures because of its potential
applications in various fields, such as communication
\cite{VanWiggeren1998}, seismology \cite{Vieira1999}, and neural
networks \cite{Hoppensteadt2000}.

The word ``synchronization" comes from a Greek word, which means
``share time". Today, in science and technology, it has come to be
considered as ``time coherence of different processes". Since the
first observation of synchronization phenomenon was made by Huygens
\cite{Huygens1672} in the 17th century, many different types of
synchronization phenomena have been found, e.g., phase
synchronization, lag synchronization, full synchronization, partial
synchronization, almost synchronization, and so on. In mathematics,
synchronization can be defined as a process wherein two (or many)
dynamical systems adjust a given property of their motion to a
common behavior as time goes to infinity, due to coupling or forcing
(see \cite{Boccaletti2002}). For example, full synchronization
requires that the difference between any two nodes converges to zero
as time goes to infinity. Therefore, it is natural to raise
following question: For given initial conditions, can the solution
be extended to the infinite interval $[0,+\infty)$?

For example, in the paper \cite{Lu2006-2}, the following coupled
systems with a delay is considered:
\begin{eqnarray}
\dot{x}^{i}(t)=f(x^{i}(t))+c\sum\limits_{j=1,j\ne
i}^{m}a_{ij}\Gamma\bigg[x^{j}(t-\tau)-x^{i}(t)\bigg],\nonumber\\
i=\onetom, \label{LCODECD1}
\end{eqnarray}
where $x^{i}(t)=[x^{i}_{1}(t),\cdots,x^{i}_{n}(t)]^{\top}\in \mathbb
R^{n}$ denotes the $n$-dimensional state variable of the $i$-th
node, $i=\onetom$; $f:\mathbb R^{n}\rightarrow \mathbb R^{n}$ is a
differential function of the intrinsic system; $c$ is the coupling
strength; $\Gamma=\mathrm{diag}\{\gamma_{1},\cdots,\gamma_{n}\}$ is
the inner connection diagonal matrix with $\gamma_{i}\ge 0$,
$i=\oneton$; $a_{ij}\ge 0$, for all $i\ne j$, is the coupling
coefficient from node $j$ to node $i$; and $\tau\ge 0$ is the
coupling delay. It is assumed that $\sum_{j=1,j\ne i}^{m}a_{ij}=1$,
$a_{ii}=-1$, for all $i=\onetom$. And the following theorem was
proved.

\begin{myprop} 
Suppose that there are a positive definite diagonal matrix
$P=\mathrm{diag}\{p_{1},\cdots,p_{n}\}$ and a diagonal matrix
$D=\mathrm{diag}\{d_{1},\cdots,d_{n}\}$, such that
\begin{eqnarray*}
(x-y)^{\top}P[f(x)-f(y)-Dx+Dy]\le -\alpha (x-y)^{\top}(x-y)
\end{eqnarray*}
holds for some $\alpha>0$, any $x,y\in \mathbb R^{n}$. Then, for
sufficiently large coupling strength $c$ and sufficiently small
delay $\tau$, the coupled system (\ref{LCODECD1}) will be globally
synchronized.
\end{myprop}

Here, a prerequisite condition in discussing synchronization is that
the solution $x^i(t)$, $i=1,\cdots,m$, can be extended to the
infinite interval $[0,+\infty)$. However, in most papers on
synchronization of coupled systems, such as
\cite{Wu1995,Lu2004,Zhou2006,Lu2006-2} and others, it is always
assumed that for each initial condition, the coupled system under
consideration has a unique solution for all time $t\geq0$ without
any theoretical justification.

In this short paper, we address this issue and propose a general
model of coupled dynamical systems, which includes previously
studied systems as special cases. We prove that under the assumption
of QUAD (Assumption (A5) in Section \ref{sec:model}), the solution
of the general model exists on $[0,+\infty)$. The assumption of QUAD
is often used when using a Lyapunov function with a quadratic form
to investigate the global synchronization (e.g., in Proposition 1,
and in \cite{Wu1995,Belykh2006,Lu2006-1}). Therefore, the theorem
proved in this paper provides a theoretical basis for the discussion
of synchronization of the coupled systems.

The rest of the paper is organized as follows: In Section
\ref{sec:model}, we propose a general model of coupled dynamical
systems. In Section \ref{sec:preliminaries}, we present some
fundamental theorems of retarded functional differential equations
with infinite delay, which are taken from \cite{Hino1991}. In
Section \ref{sec:main-result}, the main theorem is proved. We
conclude the paper in Section \ref{sec:conclusions}.

\section{Model descriptions}
\label{sec:model}

In this section, we investigate the coupled dynamical systems
described by the following retarded functional integro-differential
equations:
\begin{eqnarray}
\dot{x}^i(t) & = &
f(t,x^i(t))\nonumber\\
&   &
+\sum_{j=1}^{m}a_{ij}(t)\int_0^{\infty}g(t,x^j(t-\tau_{ij}(t)-s))
\mathrm{d}K_{ij}(s),\nonumber\\
&   &
i=1,2,\ldots,m, \label{general-model}
\end{eqnarray}
where ``$\dot{\;\;}$" represents the right-hand derivative, $m$ is
the network size, $x^i(t)\in \mathbb{R}^n$ is the state variable of
the $i$-th node, $t\in[0,+\infty)$ is a continuous time,
$f:[0,+\infty)\times\mathbb{R}^n\rightarrow\mathbb{R}^n$ describes
the dynamical behavior of each uncoupled system,
$A(t)=(a_{ij}(t))\in\mathbb{R}^{m\times m}$ is the time-varying
coupling matrix, which is determined by the topological structure of
the network,
$g:[0,+\infty)\times\mathbb{R}^n\rightarrow\mathbb{R}^n$ is the
output function, $\mathrm{d}K_{ij}(s)$ is a Lebesgue-Stieljies
measure for each $i,j=1,\ldots,m$, and satisfies
$\int_0^{\infty}|\mathrm{d}K_{ij}(s)|<+\infty$.

In addition, the following assumptions are necessary in discussion
of retarded systems:
\begin{list}
{{\upshape (A\arabic{mycounter})}\hfill} {\setlength{\topsep}{0ex}
 \setlength{\parskip}{0ex}
 \setlength{\itemsep}{0.2ex}
 \setlength{\parsep}{0ex}
 \setlength{\leftmargin}{6ex}
 \setlength{\labelwidth}{4ex}
 \setlength{\labelsep}{1ex}
 \setlength{\itemindent}{-1ex}
 \usecounter{mycounter}}
\item $f(t,u)$ is continuous, and locally Lipschitz continuous with
respect to $u$, i.e., in each compact subset $W$ of
$[0,+\infty)\times\mathbb{R}^n$, there exists a constant $l(W)>0$
such that $\|f(t,u_1)-f(t,u_2)\|\leq l(W)\|u_1-u_2\|$ for any
$(t,u_k)\in W$, $k=1,2$;

\item $A(t)=(a_{ij}(t))_{i,j=1}^m$ is continuous;

\item $g(t,u)$ is continuous, and there exists a continuous function
$\kappa(t):[0,+\infty)\rightarrow\mathbb{R}^{+}$, such that
$\|g(t,u_1)-g(t,u_2)\|\leq\kappa(t)\|u_1-u_2\|$ for any
$t\in[0,+\infty)$ and $u_1,u_2\in\mathbb{R}^n$;

\item For each $i,j=1,\ldots,m$, $\tau_{ij}(t)$ is
continuous and nonnegative;

\item There are a symmetric positive definite matrix
$P$ and a diagonal matrix
$\Delta=\mathrm{diag}\{\delta_1,\ldots,\delta_n\}$ such that
$f(t,u)\in\mathrm{QUAD}(\Delta,P)$, where $\mathrm{QUAD}(\Delta,P)$
denotes a class of continuous functions
$h(t,u):[0,+\infty)\times\mathbb{R}^n\rightarrow\mathbb{R}^n$
satisfying
\begin{eqnarray}
(u_1-u_2)^{\top}P\{[h(t,u_1)-h(t,u_2)]-\Delta[u_1-u_2]\}\nonumber\\
\leq-\epsilon(u_1-u_2)^{\top}(u_1-u_2)
\label{QUAD}
\end{eqnarray}
for some $\epsilon>0$, all $u_1,u_2\in\mathbb{R}^n$ and
$t\in[0,+\infty)$.
\end{list}
Here, $\|\cdot\|$ can be any norm in $\mathbb{R}^n$ (Without loss of
generality, in this paper we assume that $\|\cdot\|$ is 2-norm).

The model (\ref{general-model}) includes many previously studied
systems as special cases. In the following, we present several
examples.

\begin{example} 
$\mathrm{d}K_{ij}(s)=\delta(s)$, where $\delta(s)$ is the
Dirac-delta function, i.e., $\delta(0)=1$ and $\delta(s)=0$ for
$s\neq0$; $A(t)=A$ is a constant matrix with zero-sum rows and
nonnegative off-diagonal elements; $g(t,u)=\Gamma u$, where $\Gamma$
is a constant matrix; $\tau_{ij}(t)=0$ for each $i,j=1,\ldots,m$ and
all $t\geq0$. Then, (\ref{general-model}) reduces to the system with
undelayed, constant and linear coupling discussed in
\cite{Wu1995,Lu2006-1}:
\begin{eqnarray*}
\dot{x}^i(t)=f(t,x^i(t))+\sum_{j=1}^{m}a_{ij}\Gamma x^j(t), \qquad
i=1,2,\ldots,m.
\end{eqnarray*}
\end{example}

\begin{example} 
$\mathrm{d}K_{ij}(s)=\delta(s)$, where $\delta(s)$ is the
Dirac-delta function; $A(t)$ is a time-dependent matrix with
zero-sum rows and nonnegative off-diagonal elements;
$g(t,u)=\Gamma(t)u$, where $\Gamma(t)$ is a time-dependent matrix;
$\tau_{ij}(t)=0$ for each $i,j=1,\ldots,m$ and all $t\geq0$. Then,
(\ref{general-model}) reduces to the system with undelayed,
time-varying and linear coupling discussed in \cite{Wu2003,Wu2005}:
\begin{eqnarray*}
\dot{x}^i(t)=f(t,x^i(t))+\sum_{j=1}^{m}a_{ij}(t)\Gamma(t)x^j(t),\\
\qquad i=1,2,\ldots,m.
\end{eqnarray*}
\end{example}

\begin{example} 
$\mathrm{d}K_{ij}(s)=\delta(s)$, where $\delta(s)$ is the
Dirac-delta function; $f(t,u)=f(u)$, i.e., $f$ is independent of
$t$; $A(t)=A$ is a constant matrix with zero-sum rows and
nonnegative off-diagonal elements, and satisfies $a_{ii}=-c$ for
$i=1,\ldots,m$; $g(t,u)=\Gamma u$, where $\Gamma$ is a diagonal
matrix with nonnegative diagonal elements; $\tau_{ij}(t)=\tau$ for
$i\neq j$ and $\tau_{ii}(t)=0$ for $i=1,\ldots,m$. Then,
(\ref{general-model}) reduces to the system with delayed, constant
and linear coupling discussed in \cite{Lu2006-2}:
\begin{eqnarray*}
\dot{x}^i(t)=f(x^i(t))+\sum_{j=1,j\neq
i}^{m}a_{ij}\Gamma[x^j(t-\tau)-x^i(t)],\\
\qquad i=1,2,\ldots,m.
\end{eqnarray*}
\end{example}

Besides Examples 1-3, the model (\ref{general-model}) includes
coupled dynamical systems with nonlinear coupling, time-varyingly
delayed coupling, distributedly delayed coupling, etc.

\section{Preliminaries}
\label{sec:preliminaries}

In this section, we present some fundamental results of retarded
functional differential equations with infinite delay, which will be
used in the sequel.

Firstly, we introduce some notations and definitions.

Denote $BC((-\infty,a],\mathbb{R}^N)$ the family of continuous
functions $\phi$ mapping the interval $(-\infty,a]$ into
$\mathbb{R}^N$ such that
$\|\phi\|=:\sup\{\|\phi(\theta)\|:-\infty<\theta\leq a\}$ is finite.
Also, denote $C^{\infty}((-\infty,a],\mathbb{R}^N)=\{\phi\in
BC((-\infty,a],\mathbb{R}^N):\lim_{\theta\rightarrow-\infty}\phi(\theta)$
exists in $\mathbb{R}^N\}$. When $a=0$, we generally denote
$C^{\infty}=C^{\infty}((-\infty,0],\mathbb{R}^N)$. For
$\sigma\in\mathbb{R}$, $B\geq0$, $x\in
C^{\infty}((-\infty,\sigma+B],\mathbb{R}^N)$, and
$t\in[\sigma,\sigma+B]$, we define $x_t\in C^{\infty}$ as
$x_t(\theta)=x(t+\theta)$, $\theta\in(-\infty,0]$. Assume $\Omega$
is an open subset of $\mathbb{R}\times C^{\infty}$, $h:
\Omega\rightarrow\mathbb{R}^N$ is a given function, and
``$\dot{\;\;}$" represents the right-hand derivative; then, we call
\begin{eqnarray}
\dot{x}(t)=h(t,x_t) \label{RFDE}
\end{eqnarray}
a retarded functional differential equation with infinite delay on
$\Omega$.

\begin{definition} 
A function $x$ is said to be a solution of Equation (\ref{RFDE}) on
the interval $I=[\sigma,\sigma+B)$ if there are
$\sigma\in\mathbb{R}$ and $B>0$ such that $x\in
C^{\infty}((-\infty,\sigma+B),\mathbb{R}^N)$, $(t,x_t)\in\Omega$ and
$x(t)$ satisfies Equation (\ref{RFDE}) for $t\in I$. For given
$\sigma\in\mathbb{R}$, $\varphi\in C^{\infty}$, if a solution $x$ of
Equation (\ref{RFDE}) is defined on an interval $[\sigma,\sigma+B)$,
$B>0$, and satisfies $x_{\sigma}=\varphi$, then $x$ is called a
solution of Equation (\ref{RFDE}) with initial value $\varphi$ at
$\sigma$ or simply a solution through $(\sigma,\varphi)$.
\end{definition}

\begin{definition} 
Suppose $x(t)$ and $y(t)$ are solutions with the same initial
condition and satisfies Equation (\ref{RFDE}) respectively on the
intervals $I$ and $J$ whose left end points are $\sigma$. If $I$ is
properly contained in $J$ and $x(t)=y(t)$ for $t\in I$, we say $y$
is a continuation of $x$. If $x$ has no continuation, it is called a
noncontinuable solution, or a maximal solution.
\end{definition}

\begin{definition} 
We say $h(t,\phi)$ is Lipschitz in $\phi$ in a compact subset $W$ of
$\mathbb{R}\times C^{\infty}$ if there a constant $l>0$ such that,
for any $(t,\phi_k)\in W$, $k=1,2$,
\begin{eqnarray}
\|h(t,\phi_1)-h(t,\phi_2)\|\leq l\|\phi_1-\phi_2\|.\label{Lip}
\end{eqnarray}
\end{definition}

The following three lemmas on existence, uniqueness, and
continuation of the solution of Equation (\ref{RFDE}),  are used in
the proof of the main theorem in the next section. The details can
be found in \cite{Hino1991},

\begin{lemma} 
(Existence) Suppose $\Omega$ is an open subset in
$\mathbb{R}\times C^{\infty}$ and $h:\Omega\rightarrow\mathbb{R}^N$
is continuous. Then, for any $(\sigma,\varphi)\in\Omega$, there
exists a solution of Equation (\ref{RFDE}) through
$(\sigma,\varphi)$.
\end{lemma}

\begin{lemma} 
(Uniqueness) Suppose $\Omega$ is an open subset in
$\mathbb{R}\times C^{\infty}$ and $h(t,\phi)$ is Lipschitz in $\phi$
in each compact subset of $\Omega$. Then, for any
$(\sigma,\varphi)\in\Omega$, there exists at most one noncontinuable
solution of Equation (\ref{RFDE}) through $(\sigma,\varphi)$.
\end{lemma}

\begin{lemma} 
(Continuation) Suppose $\Omega$ is an open subset in
$\mathbb{R}\times C^{\infty}$, $h:\Omega\rightarrow\mathbb{R}^N$ is
continuous, and $x$ is a noncontinuable solution of Equation
(\ref{RFDE}) defined on $I=[\sigma,\sigma+B)$. Then, for every
compact subset $W$ of $\Omega$, there is a $t_W$ in $I$ such that
$(t,x_t)\not\in W$ for all $t\in(t_W,\sigma+B)$.
\end{lemma}

\section{Main result}
\label{sec:main-result}

In this section, we prove the following theorem.

\begin{theorem} 
Suppose that Assumptions (A1)-(A5) hold. Then, for any
$\varphi(\theta)=[\varphi^1(\theta)^{\top},\ldots,\varphi^m(\theta)^{\top}]^{\top}$
with $\varphi^i(\theta)\in C^{\infty}((-\infty,0],\mathbb{R}^n)$,
there is a unique noncontinuable solution
$x(t)=[x^1(t)^{\top},\ldots,x^m(t)^{\top}]^{\top}$ of Equation
(\ref{general-model}) through $(0,\varphi)$. Moreover, the interval
of existence of the solution $x$ is $[0,+\infty)$.
\end{theorem}

{\it Proof :} By Assumptions (A1)-(A4) and Lemmas 1-2, it is clear
that for the integro-diffential system (\ref{general-model}), there
exists a unique noncontinuable solution $x(t)$. In the following, we
will prove that the interval of existence of the solution $x(t)$ is
$[0,+\infty)$.

We employ ``proof by contradiction", and suppose that the interval
of existence of the noncontinuable solution $x(t)$ is $[0,b)$, where
$b$ is a positive constant.

Firstly, by Assumptions (A1)-(A4), we can find positive constants
$\alpha$, $\beta$ and $\gamma$ such that
\begin{eqnarray*}
&&\|g(t,u_1)-g(t,u_2)\|\leq\alpha\|u_1-u_2\|
\end{eqnarray*}
holds for all $u_1,u_2\in\mathbb{R}^n$ and $t\in[0,b)$, and
\begin{eqnarray*}
&&|a_{ij}(t)|\leq\beta,\\
&&\big\|f(t,x^i(0))+\sum_{j=1}^m
a_{ij}(t)g(t,x^j(0))\int_0^{\infty}\mathrm{d}K_{ij}(s)\big\|\leq\gamma.
\end{eqnarray*}
hold for all $i,j=1,\ldots,m$ and $t\in[0,b)$.

Now, we will show how the assumption of QUAD (Assumption (A5)) plays
an important role in the proof.

Since $f(t,u)\in\mathrm{QUAD}(\Delta,P)$ (Assumption (A5)), it is
clear that there is a constant $\delta>0$ such that for all
$u_1,u_2\in\mathbb{R}^n$ and $t\geq0$,
\begin{eqnarray*}
(u_1-u_2)^{\top}P[f(t,u_1)-f(t,u_2)]\leq\delta(u_1-u_2)^{\top}(u_1-u_2).
\end{eqnarray*}
Denote
\begin{eqnarray*}
\eta=\frac{2\delta+2\alpha\beta\big\|P\big\|K}{\lambda^P_{\min}}
+\frac{2m\gamma\big\|P\big\|}{\sqrt{\lambda^P_{\min}}}>0\;,
\end{eqnarray*}
where $m$ is the number of the nodes, $\|P\|$ is the 2-norm of the
matrix $P$, $\lambda^P_{\min}$ is the minimum eigenvalue of the
matrix $P$, and
$K=\sum_{i=1}^m\sum_{j=1}^m\int_0^{\infty}|\mathrm{d}K_{ij}(s)|$.

Since the matrix $P$ is symmetric positive definite, we can define a
norm in $\mathbb{R}^{nm}$:
\begin{eqnarray*}
\|x(t)\|_P=\Big(\sum_{i=1}^m
x^i(t)^{\top}Px^i(t)\Big)^{\frac{1}{2}};
\end{eqnarray*}
and two nonnegative functions:
\begin{eqnarray*}
&&V(t)=\frac{1}{2}\|x(t)-x(0)\|_P^2,\\
&&M(t)=\max\Big[\frac{1}{2}\,,\; \sup_{-\infty<s\leq
t}\frac{1}{2}\|x(s)-x(0)\|_P^2\Big],\\
&&t\in[0,b).
\end{eqnarray*}
Clearly, $V(t)\leq M(t)$.We claim that $M(t)\leq M(0)e^{\eta t}$ for
all $t\in[0,b)$.

In fact, at any $t_0\in[0,b)$, there are two possible cases:

\noindent{\bf Case 1:} $V(t_0)<M(t_0)$. In this case, by the
continuity of $\|x(t)-x(0)\|_P^2$, $M(t)$ is non-increasing at
$t_0$.

\noindent{\bf Case 2:} $V(t_0)=M(t_0)$.

Calculating the right-hand derivative of $V$ with respect to time
along the trajectories of (\ref{general-model}), one has
\allowdisplaybreaks[3]
\begin{eqnarray*}
&      &
\hspace{-2em}\dot{V}(t_0)=
\sum_{i=1}^m(x^i(t_0)-x^i(0))^{\top}P\dot{x^i}(t_0)\allowdisplaybreaks[3]\\
&  =   &
\sum_{i=1}^m(x^i(t_0)-x^i(0))^{\top}P\bigg[f(t_0,x^i(t_0))\\
&      &
+\sum_{j=1}^ma_{ij}(t_0)
\int_0^{\infty}g(t_0,x^j(t_0-\tau_{ij}(t_0)-s))\mathrm{d}K_{ij}(s)\bigg]\allowdisplaybreaks[3]\\
&  =   &
\sum_{i=1}^m(x^i(t_0)-x^i(0))^{\top}P\bigg\{\big[f(t_0,x^i(t_0))-f(t_0,x^i(0))\big]\\
&      &
+\sum_{j=1}^m a_{ij}(t_0)\int_0^{\infty}\big[g(t_0,x^j(t_0-\tau_{ij}(t_0)-s))\\
&      &
-g(t_0,x^j(0))\big]\mathrm{d}K_{ij}(s)+\Big[f(t_0,x^i(0))\\
&      &
+\sum_{j=1}^ma_{ij}(t_0)g(t_0,x^j(0))\int_0^{\infty}\mathrm{d}K_{ij}(s)\Big]\bigg\}\allowdisplaybreaks[3]\\
& \leq &
\delta\sum_{i=1}^m(x^i(t_0)-x^i(0))^{\top}(x^i(t_0)-x^i(0))\\
&      &
+\sum_{i=1}^m\sum_{j=1}^m\big|a_{ij}(t_0)\big|\big\|P\big\|\big\|x^i(t_0)-x^i(0)\big\|\\
&      &
\times\int_0^{\infty}\big\|g(t_0,x^j(t_0-\tau_{ij}(t_0)-s))\\
&      &
-g(t_0,x^j(0))\big\|\big|\mathrm{d}K_{ij}(s)\big|\\
&      &
+\sum_{i=1}^m\big\|P\big\|\big\|x^i(t_0)-x^i(0)\big\|\big\|f(t_0,x^i(0))\\
&      &
+\sum_{j=1}^ma_{ij}(t_0)g(t_0,x^j(0))\int_0^{\infty}\mathrm{d}K_{ij}(s)\big\|\allowdisplaybreaks[3]\\
& \leq &
\delta\big\|x(t_0)-x(0)\big\|^2
+\alpha\beta\big\|P\big\|\sum_{i=1}^m\sum_{j=1}^m\big\|x^i(t_0)-x^i(0)\big\|\\
&      &
\times\int_0^{\infty}\big\|x^j(t_0-\tau_{ij}(t_0)-s)-x^j(0)\big\|\big|\mathrm{d}K_{ij}(s)\big|\\
&      &
+\gamma\big\|P\big\|\sum_{i=1}^m\big\|x^i(t_0)-x^i(0)\big\|\allowdisplaybreaks[3]\\
& \leq &
\delta\big\|x(t_0)-x(0)\big\|^2
+\alpha\beta\big\|P\big\|\sum_{i=1}^m\sum_{j=1}^m\big\|x(t_0)-x(0)\big\|\\
&      &
\times\int_0^{\infty}\big\|x(t_0-\tau_{ij}(t_0)-s)-x(0)\big\|\big|\mathrm{d}K_{ij}(s)\big|\\
&      &
+m\gamma\big\|P\big\|\big\|x(t_0)-x(0)\big\|\allowdisplaybreaks[3]\\
& \leq &
\frac{\delta}{\lambda^P_{\min}}\big\|x(t_0)-x(0)\big\|_P^2\\
&      &
+\frac{\alpha\beta\big\|P\big\|}{\lambda^P_{\min}}\sum_{i=1}^m\sum_{j=1}^m\big\|x(t_0)-x(0)\big\|_P\\
&      &
\times\int_0^{\infty}\big\|x(t_0-\tau_{ij}(t_0)-s)-x(0)\big\|_P\big|\mathrm{d}K_{ij}(s)\big|\\
&      &
+\frac{m\gamma\big\|P\big\|}{\sqrt{\lambda^P_{\min}}}\big\|x(t_0)-x(0)\big\|_P\cdot1\allowdisplaybreaks[3]\\
& \leq &
\frac{\delta}{\lambda^P_{\min}}2M(t_0)\\
&      &
+\frac{\alpha\beta\big\|P\big\|}{\lambda^P_{\min}}\sqrt{2M(t_0)}\sqrt{2M(t_0)}
\sum_{i=1}^{m}\sum_{j=1}^{m}\int_0^{\infty}\big|\mathrm{d}K_{ij}(s)\big|\\
&      &
+\frac{m\gamma\big\|P\big\|}{\sqrt{\lambda^P_{\min}}}\sqrt{2M(t_0)}\sqrt{2M(t_0)}\allowdisplaybreaks[3]\\
&  =   &
\bigg\{\frac{2\delta+2\alpha\beta\big\|P\big\|K}{\lambda^P_{\min}}
+\frac{2m\gamma\big\|P\big\|}{\sqrt{\lambda^P_{\min}}}\bigg\}M(t_0)\\
&  =   & \eta V(t_0)
\end{eqnarray*}

In summary, we conclude that $M(t)\leq M(0)e^{\eta t}$ for all
$t\in[0,b)$, which implies $V(t)\leq M(0)e^{\eta t}$ and
\begin{eqnarray}
\|x(t)\| & \leq & \frac{1}{\sqrt{\lambda^P_{\min}}}\|x(t)\|_P\nonumber\\
& \leq &
\frac{1}{\sqrt{\lambda^P_{\min}}}\big(\|x(0)\|_P+\|x(t)-x(0)\|_P\big)\nonumber\\
&  =   &
\frac{1}{\sqrt{\lambda^P_{\min}}}\big(\|x(0)\|_P+\sqrt{2V(t)}\;\big)\nonumber\\
& \leq &
\frac{1}{\sqrt{\lambda^P_{\min}}}\big(\|x(0)\|_P+\sqrt{2M(0)e^{\eta
t}}\;\big)\nonumber\\
& \leq &
\frac{1}{\sqrt{\lambda^P_{\min}}}\big(\|x(0)\|_P+\sqrt{2M(0)e^{\eta
b}}\;\big)\label{boundary}
\end{eqnarray}
for all $t\in[0,b)$.

Now, pick a compact set
\begin{eqnarray*}
W = \bigg\{(t,\psi)\in\mathbb{R}\times C^{\infty}((-\infty,0],
\mathbb{R}^{nm})\Big|0\leq t\leq b,\mbox{ and }\\
\|\psi\|\leq\max\Big[\frac{1}{\sqrt{\lambda^P_{\min}}}\big(\|x(0)\|_P
+\sqrt{2M(0)e^{\eta b}}\;\big),\;\|\varphi\|\Big]\bigg\},
\end{eqnarray*}
where $\varphi$ is the initial value. By the inequality
(\ref{boundary}), we conclude that $(t,x_t)\in W$ for all
$t\in[0,b)$, which contradicts Lemma 3.

Therefore, the interval of existence of the noncontinuable solution
$x$ is $[0,+\infty)$. Theorem is proved completely.

\section{Conclusions}
\label{sec:conclusions}

In this paper, we propose a general model of coupled dynamical
systems, which includes previously studied systems as special cases,
and prove that under the assumption of QUAD, the solution of the
general model exists on $[0,+\infty)$.

\end{document}